\documentclass[12pt]{article}

\usepackage{amsmath,amssymb,pb-diagram}

\begin{document}

\pagestyle{myheadings} \markright{PRIME GEODESIC THEOREM...}

\title{A prime geodesic theorem for higher rank II: singular geodesics}
\author{Anton Deitmar}

\date{}
\maketitle

\tableofcontents

\def \1{{\bf 1}}
\def \a{{{\mathfrak a}}}
\def \ad{{\rm ad}}
\def \al{\alpha}
\def \ar{{\alpha_r}}
\def \A{{\mathbb A}}
\def \Ad{{\rm Ad}}
\def \Aut{{\rm Aut}}
\def \b{{{\mathfrak b}}}
\def \bs{\backslash}
\def \B{{\rm B}}
\def \c{{\mathfrak c}}
\def \cent{{\rm cent}}
\def \C{{\mathbb C}}
\def \CA{{\cal A}}
\def \CB{{\cal B}}
\def \CC{{\cal C}}
\def \CD{{\cal D}}
\def \CE{{\cal E}}
\def \CF{{\cal F}}
\def \CG{{\cal G}}
\def \CH{{\cal H}}
\def \CHC{{\cal HC}}
\def \CL{{\cal L}}
\def \CM{{\cal M}}
\def \CN{{\cal N}}
\def \CP{{\cal P}}
\def \CQ{{\cal Q}}
\def \CO{{\cal O}}
\def \CS{{\cal S}}
\def \CT{{\cal T}}
\def \CV{{\cal V}}
\def \CW{{\cal W}}
\def \d{{\mathfrak d}}
\def \df{\ \begin{array}{c} _{\rm def}\\ ^{\displaystyle =}\end{array}\ }
\def \det{{\rm det}}
\def \diag{{\rm diag}}
\def \dist{{\rm dist}}
\def \End{{\rm End}}
\def \eqn{\begin{eqnarray*}}
\def \endeqn{\end{eqnarray*}}
\def \eps{\varepsilon}
\def \F{{\mathbb F}}
\def \Fx{{\mathfrak x}}
\def \FX{{\mathfrak X}}
\def \g{{{\mathfrak g}}}
\def \ga{\gamma}
\def \Ga{\Gamma}
\def \Gal{{\rm Gal}}
\def \h{{{\mathfrak h}}}
\def \Hom{{\rm Hom}}
\def \im{{\rm im}}
\def \Im{{\rm Im}}
\def \ind{{\rm ind}}
\def \k{{{\mathfrak k}}}
\def \K{{\cal K}}
\def \l{{\mathfrak l}}
\def \la{\lambda}
\def \lap{\triangle}
\def \li{{\rm li}}
\def \La{\Lambda}
\def \Lie{{\rm Lie}}
\def \m{{{\mathfrak m}}}
\def \mod{{\rm mod}}
\def \n{{{\mathfrak n}}}
\def \name{\bf}
\def \Mat{{\rm Mat}}
\def \N{\mathbb N}
\def \o{{\mathfrak o}}
\def \ord{{\rm ord}}
\def \O{{\cal O}}
\def \p{{{\mathfrak p}}}
\def \ph{\varphi}
\def \prf{\noindent{\bf Proof: }}
\def \Per{{\rm Per}}
\def \q{{\mathfrak q}}
\def \qed{\ifmmode\eqno $\square$\else\noproof\vskip 12pt plus 3pt minus 9pt \fi}
 \def\noproof{{\unskip\nobreak\hfill\penalty50\hskip2em\hbox{}%
     \nobreak\hfill $\square$\parfillskip=0pt%
     \finalhyphendemerits=0\par}}
\def \Q{\mathbb Q}
\def \res{{\rm res}}
\def \R{{\mathbb R}}
\def \Re{{\rm Re \hspace{1pt}}}
\def \r{{\mathfrak r}}
\def \ra{\rightarrow}
\def \rank{{\rm rank}}
\def \setminus{\begin{picture}(18,10)\put(4,6)
                {\line(2,-1){10}}\end{picture}}
\def \SL{{\rm SL}}
\def \SO{{\rm SO}}
\def \supp{{\rm supp}}
\def \Spin{{\rm Spin}}
\def \t{{{\mathfrak t}}}
\def \T{{\mathbb T}}
\def \tr{{\hspace{1pt}\rm tr\hspace{2pt}}}
\def \vol{{\rm vol}}
\def \z{\zeta}
\def \Z{\mathbb Z}
\def \={\ =\ }
\def \({\left(}
\def \){\right)}

\newcommand{\frack}[2]{\genfrac{}{}{0pt}{}{#1}{#2}}
\newcommand{\rez}[1]{\frac{1}{#1}}
\newcommand{\der}[1]{\frac{\partial}{\partial #1}}
\renewcommand{\binom}[2]{\left( \begin{array}{c}#1\\#2\end{array}\right)}
\newcommand{\norm}[1]{\parallel #1 \parallel}
\renewcommand{\matrix}[4]{\left(\begin{array}{cc}#1 & #2 \\ #3 & #4 \end{array}\right)}
\renewcommand{\sp}[2]{\langle #1,#2\rangle}
\renewcommand{\labelenumi}{(\alph{enumi})}

\newtheorem{theorem}{Theorem}[section]
\newtheorem{conjecture}[theorem]{Conjecture}
\newtheorem{lemma}[theorem]{Lemma}
\newtheorem{corollary}[theorem]{Corollary}
\newtheorem{proposition}[theorem]{Proposition}

\newpage

\begin{center} {\bf Introduction} \end{center}

The prime geodesic
theorem gives a growth asymptotic for the number of closed geodesics counted by their lengths
\cite{Gangolli,Hejhal,Knieper,Koyama,LuoSarnak,Margulis,Pollicott,Zelditch}.
It has hitherto only
been proven for manifolds of strictly negative curvature. 
For manifolds containing higher dimensional flats it is not a priori clear what a
prime geodesic theorem might look like. 
In the paper \cite{primgeo} the author has given such a theorem for locally symmetric spaces of arbitrary rank, i.e., they may contain
higher dimensional flats. The regular geodesics in such a space give points in a
higher dimensional Weyl cone, and the prime geodesic theorem describes the distribution of these points. 
In the current paper we turn to the remaining, i.e., singular geodesics.
As already mentioned in \cite{primgeo}, there are serious obstacles to giving an asymptotical formula in general, but if one imposes extra regularity conditions on the space, then these obstacles disappear and one can derive an asymptotical formula for singular geodesics.

We describe the main result of the paper. One of the
various equivalent formulations of the prime geodesic theorem for locally symmetric spaces of rank one is
the following. Let $\bar X$ be a compact locally symmetric
space with universal covering of rank one. For $T>0$ let
$$
\psi(T)\=\sum_{c\,:\, e^{l(c)}\le T} l(c_0).
$$
Here the sum runs over all closed geodesics $c$ such that
$e^{l(c)}\le T$, where $l(c)$ is the length of the
geodesic $c$, and $c_0$ is the prime geodesic underlying
$c$. Then, under a suitable scaling of the metric, as $T\ra\infty$,
$$
\psi(T)\ \sim\ T.
$$

We now replace the space $\bar X$ by an arbitrary compact locally symmetric
space which is a quotient of a globally symmetric space $X=G/K$ where $G$ is a semisimple
Lie group of split-rank
$r$ and
$K$ a maximal compact subgroup.  
So the space under consideration is $\Ga\bs X-\Ga\bs G/K$, where $\Ga\subset G$ is a torsion-free discrete subgroup.
The extra regularity condition one has to put on this space is that $\Ga$ be a \emph{regular} group (see section 1).
A closed geodesic $c$ gives rise to
a point $a_c$ in the closure of the negative Weyl chamber $A_0^-$ of a maximal split torus $A_0$.  
We pick a wall $A^-$ of $A_0^-$ which might be equal to $A_0^-$ or of smaller dimension.
We consider all geodesics $c$ that give points $a_c$ in $A^-$.
Let $r$ be the dimension of $A^-$.
For
$T_1,\dots T_r>0$ let
$$
\psi(T_1,\dots T_r)\=\sum_{c\,:\, a_{c,j}\le T_j} \la_c,
$$
where $\la_c$ is the volume of the unique maximal flat $c$ lies in and $a_{c,j}$ are the coordinates of $a_c$ with
respect to a canonical coordinate system on $A^-$ given by the roots. The sum runs over all
closed geodesics $c$ with $a_c\in A^-$ modulo homotopy. The main result of this paper is that,
as $T_j$  tends to infinity for every $j$,
$$
\psi(T_1,\dots T_r)\ \sim\ T_1\cdots T_r.
$$
The  proof is based on a Lefschetz formula similar to the one \cite{primgeo}, but at various places one has to argue in a fashion different to the previous case.

The restriction to regular spaces is a strong one, but fortunately the most important application which is an asymptotic formula for units in orders, can be derived in this context if the degree of the number field generated by the order is a prime.

\section{The Lefschetz formula}\label{sec1}
In this section we give a Lefschetz formula for regular locally symmetric spaces.  
Let $G$ be a connected semisimple Lie group with finite
center and choose a maximal compact subgroup $K$ with
Cartan involution $\theta$, i.e., $K$ is the group of
fixed points of $\theta$. Let $P$ be a cuspidal parabolic subgroup with Langlands decomposition
$P=MAN$. 
Cuspidality here means that the group $M$ admits a compact Cartan subgroup.
Modulo conjugation we can assume that $A$ and $M$ are stable under
$\theta$. 
The centralizer of $A$
is $AM$. Let $W(A,G)$ be the \emph{Weyl group} of $A$, i.e. $W(A,G)$ is the quotient of the
normalizer of $A$ by the centralizer. 
This is a finite group acting on $A$. 

We have to fix Haar measures. We use the normalization of
Harish-Chandra \cite{HC-HA1}. Note that in this normalization of Haar measures the
compact group $K$ has volume one.

We write $\g_\R, \k_\R,\a_\R,\m_\R,\n_\R$ for the real Lie
algebras of $G,K,A,M,N$ and $\g,\k,\a,\m,\n$ for their
complexifications. $U(\g)$ is the universal enveloping
algebra of $\g$. 
This algebra is isomorphic to the algebra
of all left invariant differential operators on $G$ with
complex coefficients. 
Pick a compact Cartan subgroup $T$ of $M$ and let $\t$ be its complexified Lie algebra. 
Then $\h=\a\oplus\t$ is a
Cartan subalgebra of $\g$. Let $W(\h,\g)$ be the corresponding absolute Weyl group.

Let $\a^*$ denote the dual space of the complex vector
space $\a$. Let $\a_\R^*$ be the real dual of $\a_\R$. We identify $\a_\R^*$ with the real vector space of all $\la\in\a^*$ that map $\a_\R$ to $\R$. Let $\Phi\subset\a^*$ be the set of all roots
of the pair $(\a,\g)$ and let $\Phi^+$ be the subset of
positive roots with respect to $P$. Let $\Delta\subset
\Phi^+$ be the set of simple roots. Then $\Delta$ is a basis of $\a^*$. The open {\it
negative Weyl chamber} $\a_\R^-\subset\a_\R$ is the cone of all
$X\in\a_\R$ with $\alpha(X)<0$ for every $\alpha\in\Delta$. Let $\overline{\a_\R^-}$ be the
closure of $\a_\R^-$. 

The bilinear form $B$ is indefinite on $\g_\R$, but the form
$$
\sp{X}{Y}\df -B(X,\theta(Y))
$$
is positive definite, ie an inner product on $\g_\R$. We extend it to an inner product on
the complexification $\g$. Let $\norm X=\sqrt{\sp XX}$ be the corresponding norm.
The form $B$, being nondegenerate, identifies $\g$ to its dual space $\g^*$. In this way we
also define an inner product $\sp ..$ and the corresponding norm on $\g^*$.
Furthermore, if $V\subset\g$ is any subspace on which $B$ is nondegenerate, then $B$ gives
an identification of $V^*$ with $V$ and so one gets an inner product and a norm on
$V^*$. This in particular applies to $V=\h$, a Cartan subalgebra of $\g$, which is defines
over
$\R$.

Let $\Ga\subset G$ be a discrete, cocompact, torsion-free subgroup. 
We are interested in the closed geodesics on the locally symmetric space $X_\Ga=\Ga\bs
X=\Ga\bs G/K$. Every such geodesic 
$c$ lifts to a $\Ga$-orbit of geodesics on $X$ and gives a $\Ga$-conjugacy class $[\ga_c]$ of
elements closing the particular geodesics. This induces a bijection between the set of all
homotopy classes of closed geodesics in $X_\Ga$ and the set of all non-trivial conjugacy
classes in $\Ga$ (see \cite{DKV}).

Let $C$ be an arbitrary Cartan subgroup of $G$.
The \emph{regular elements} of $C$ are
$$
C^{\rm reg}\df \{ x\in C : G_x=C\},
$$
where $G_x$ denotes the centralizer of $x$ in $G$.
Then $G^{\rm reg}$ is by definition the union of all $C^{\rm reg}$ over all Cartan subgroups.
This is an open dense set in $G$.
The group $\Ga$ is called \emph{regular} if 
$$
\Ga\setminus \{1\}\ \subset\ G^{\rm reg}.
$$
We will from now on assume that $\Ga$ is regular.

\subsection{Lefschetz numbers}
Let $\Ga$ be a torsion-free regular subgroup of $G$.
Let $\CE_P(\Ga)$ denote the set of all $\Ga$-conjugacy classes $[\ga]$ such that $\ga$ is $G$-conjugate to an element $a_\ga t_\ga$ of $A^- T$.
Let $\CE_P^0(\Ga)$ be the set of all $\Ga$-conjugacy classes $[\ga]$ such that $\ga$ is $G$-conjugate to an element $a_\ga t_\ga$ of $A^-\tilde T$, where $\tilde T$ is the intersection of $T$ with the connected component $M^0$ of the unit in $M$.
Then $\CE_P^0(\Ga)$ is a subset of $\CE_P(\Ga)$.
Let $n=\# (T/\tilde T)\in\N$, then for every $[\ga]\in\CE_P(\Ga)$ we have $[\ga^n]\in\CE_P^0(\Ga)$.

Let $[\ga]\in\CE_P(\Ga)$.
There is a closed geodesic $c$ in the Riemannian manifold $\Ga\bs G/K$ which gets closed by
$\ga$. This means that there is a lift $\tilde c$ to the universal covering $G/K$ which is
preserved by $\ga$ and $\ga$ acts on $\tilde c$ by a translation.
The closed geodesic $c$ is not unique in general. 
Since $\Ga$ is regular, there is a unique maximal flat $F_c$ containing $c$.
By \emph{maximal flat} we here mean a flat, totally geodesic submanifold which is maximal with these properties with respect to inclusion.
Note that other authors sometimes insist that a maximal flat should be of maximal dimension which we do not.
Let $\la_\ga$ be the volume of that flat,
$$
\la_\ga\df \vol(F_c).
$$ 
As the notation indicates, this number only depends on $\ga$ and not on $c$.

Let $\n$ denote the complexified Lie algebra of $N$. For
any $\n$-module $V$ let $H_q(\n,V)$ and $H^q(\n,V)$ for $q=0,\dots,\dim\n$
be the Lie algebra homology and cohomology \cite{BorWall}.  Let $\hat G$ denote the unitary dual of $G$, i.e., the set of isomorphism classes of irreducible unitary representations of $G$.  For $\pi\in\hat G$ let $\pi_K$ be the $(\g,K)$-module of $K$-finite vectors. If
$\pi\in\hat G$, then $H_q(\n,\pi_K)$ and $H^q(\n,\pi_K)$ are
admissible $(\a\oplus\m,M)$-modules of finite length \cite{HeSch}. 

Note that $AM$ acts on the Lie algebra $\n$ of $N$ by the
adjoint representation. Let $[\ga]\in\CE_P(\Ga)$. Since
$a_\ga\in A^-$ it follows that every eigenvalue of $a_\ga
t_\ga$ on $\n$ is of absolute value $<1$. Therefore
$\det(1-a_\ga t_\ga | \n)\ne 0$.

For $[\ga]\in\CE_P(\Ga)$ let
$$
\ind(\ga)\=\frac{\la_\ga}{\det(1-a_\ga t_\ga \mid \n)}\ >\ 0,
$$
where $r=\dim A$.
Since $\Ga$ is cocompact, the unitary $G$-representation on $L^2(\Ga\bs G)$ splits
discretely with finite multiplicities
$$
L^2(\Ga\bs G)\= \bigoplus_{\pi\in\hat G} N_\Ga(\pi)\pi,
$$
where $N_\Ga(\pi)$ is a non-negative integer and $\hat G$ is the unitary dual of $G$.  A {\it quasi-character} of $A$ is a
continuous group homomorphism to $\C^\times$. Via
differentiation the set of quasi-characters can be
identified with the dual space $\a^*$. For $\la\in\a^*$ we write $a\mapsto a^\la$ for the
corresponding quasicharacter on $A$. We denote by
$\rho\in\a^*$ the modular shift with respect to $P$, i.e., for $a\in A$ we have
$\det(a|\n)=a^{2\rho}$.

For a complex vector space $V$ on which $A$ acts linearly and $\la\in\a^*$ let $(V)_\la$
denote the generalized $(\la+\rho)$-eigenspace, i.e.,
$$
(V)_\la\= \{ v\in V\mid (a-a^{\la+\rho}Id)^n v=0\ \ {\rm for\ some\ } n\in\N\}.
$$
Since $H^p(\n,\pi_K)$ is of finite length as $(\a\oplus\m,K_M)$-module, one has
$$
H^p(\n,\pi_K) \= \bigoplus_{\nu\in\a^*} H^p(\n,\pi_k)_\nu.
$$
Let $T$ be a compact Cartan subgroup of $M$ and let $\t$ be its complex Lie algebra. Then $AT$ is a
Cartan subgroup of $G$. 
Let $K_M= M\cap K$. This is a maximal compact subgroup of $M$.
Let $\Lambda_\pi\in(\a\oplus\t)^*$ be a representative of the
infinitesimal character of $\pi$. By Corollary 3.32 of \cite{HeSch} it follows,
$$
H_p(\n,\pi_K)\=\bigoplus_{\nu=w\Lambda_{\pi}|_\a}H_p(\n,\pi_K)_\nu,
$$
where $w$ ranges over $W(\g,\h)$.

\begin{lemma}\label{1.1}
For $0\le p\le d=\dim(\n)$ we have
$$
 H_p(\n,\pi_K)\ \cong\ H^{d-p}(\n,\pi_K)\otimes\det(\n),
$$
where the determinant of a
finite dimensional space is the top exterior power. So $\det(\n)$ is a one dimensional
$AM$-module on which $AM$ acts via the quasi-character $am\mapsto \det(am|\n)=a^{2\rho}$.
This in particular implies
$$
H^p(\n,\pi_K)\=\bigoplus_{\nu= w\Lambda_\pi |_\a}H^p(\n,\pi_K)_{\nu-2\rho}.
$$
\end{lemma}
\prf The first part follows straight from the definition of Lie algebra cohomology. The
second part by Corollary 3.32 of \cite{HeSch}.
\qed

Let $\m=\k_M\oplus\p_M$ be the Cartan decomposition of the complex Lie algebra $\m$ of $M$ with respect to $K_M$.
For $\la\in\a^*$
 and $\pi\in\hat G$ let
$$
m_\la(\pi)\= \sum_{q=0}^{\dim \n} \sum_{p=0}^{\dim\p_M} 
(-1)^{q+\dim\n} \dim\left(
H^q(\n,\pi_K)_\la\otimes\bigwedge^p\p_M\right)^{K_M},
$$
where the superscript $K_M$ indicates the subspace of $K_M$-invariants. Then 
$m_\la(\pi)$ is an integer and by the above, the set of $\la$ for which
$m_\la(\pi)\ne 0$ for a given $\pi$ has at most $|W(\g,\h)|$ many elements. 

Likewise define
$$
m_\la^0(\pi)\= \sum_{q=0}^{\dim \n} \sum_{p=0}^{\dim\p_M} 
(-1)^{q+\dim\n} \dim\left(
H^q(\n,\pi_K)_\la\otimes\bigwedge^p\p_M\right)^{K_M^0},
$$
where $K_M^0$ is the connected component of the unit in $K_M$, or $K_M^0=K_M\cap M^0$.

For  $\mu\in\a^{*}$ and $j\in\N$ let $\CC^{j,\mu,-}(A)$ denote the space of functions
$\ph$ on $A$ which \nopagebreak
\begin{itemize}
  \item are $j$-times continuously differentiable on $A$,
  \item are zero outside $A^-$,
  \item are such that $a^{-\mu} D\ph(a)$ is bounded on $A$  for every invariant differential operator $D$ on $A$ of
degree $\le j$.
\end{itemize}

For every invariant differential operator $D$ of degree $\le j$ let $N_D(\ph)=\sup_{a\in A}\left| a^{-\mu}\, D\ph(a)\right|$.
Then $N_D$ is a seminorm.
Let $D_1,\dots,D_n$ be a basis of the space of invariant differential operators of degree $\le j$, then
$N(\ph)=\sum_{j=1}^n N_{D_j}(\ph)$ is a norm that makes $\CC^{j,\mu,-}(A)$ into a Banach space.
A different choice of basis will give an equivalent norm.

Let $f_{EP}\in C_c^\infty(M)$ denote an Euler-Poincar\'e function on $M$.
This means that for every irreducible unitary representation $\eta$ of $M$ one has
$$
\tr\eta(f_{EP})\=\sum_{q=0}^{\dim\p_M}(-1)^q\dim\(\eta\otimes\bigwedge^p\p_M\)^{K_M}.
$$
Euler-Poincar\'e functions have the property that their orbital integrals filter out elliptic elements, i.e., for $x\in M$ a regular element one has
$$
\CO_x^M(f_{EP})\=\int_{M/M_x}f_{EP}(yxy^{-1})\, dy
$$
equals $1$ if $x$ is elliptic and zero otherwise.

\begin{theorem}\label{lefschetz} (Lefschetz Formula)\\
Assume $\Ga$ is neat.
There exists $j\in\N$ and $\mu\in\a^{*}$ such that
for any $\ph\in \CC^{j,\mu,-}(A)$ we have
$$
\sum_{\pi\in\hat G} N_\Ga(\pi)\sum_{\la\in\a^*} m_\la(\pi) \int_{A^-}\ph(a)
a^{\la+\rho} da\=\sum_{[\ga]\in\CE_P(\Ga)}\ind(\ga)\,\ph(a_\ga),
$$
where all sums and integrals converge absolutely. The inner sum on the left is always finite,
more precisely it has length $\le |W(\h,\g)|$. The left hand side is called the
\emph{global side} and the other the \emph{local side} of the Lefschetz Formula. 
Both sides of the formula give a continuous linear functional on the Banach space $\CC^{j,\mu,-}(A)$.

We also obtain a \emph{weak Lefschetz formula} by replacing $m_\la(\pi)$ with $m_\la^0(\pi)$ and $\CE_P(\Ga)$ with $\CE_P^0(\Ga)$.
\end{theorem}

\prf
The proof is in section 4 of \cite{hr} or, in a special case, in \cite{primgeo}.
The proof of the weak version is a variant of that proof where one replaces the Euler-Poincar\'e function of $M$ with the Euler-Poincar\'e function of the connected component $M^0$.
\qed

\section{The Dirichlet series}\label{sec2}
We keep assuming that the torsion-free group $\Ga$ also is regular.
Let $r=\dim A$ and for $k=1,\dots ,r$ let $\al_k$ be a positive real multiple of a
simple root of $(A,P)$ such that the
modular shift $\rho$ satisfies
$$
2\rho\= \al_1+\dots+\al_r.
$$
This defines $\al_1,\dots,\al_r$ uniquely up to order.
For $a\in A$ and $k=1,\dots r$ let $l_k(a)=|\al_k(\log a)|$ and $l(a)=l_1(a)\cdots l_r(a)$.
For $s=(s_1,\dots,s_r)\in\C^r$ and $j\in\N$ define
$$
L^j(s)\=\sum_{[\ga]\in\CE_P(\Ga)} \ind(\ga)\, 
l(a_\ga)^{j+1}\,a_\ga^{s\cdot
\al},
$$
where $s\cdot \al=s_1\al_1+\dots +s_r\al_r$. We will show that this series converges if
$\Re(s_k)>1$ for $k=1,\dots,r$. 
Likewise we define
$$
L^{0,j}(s)\=\sum_{[\ga]\in\CE_P^0(\Ga)} \ind(\ga)\, 
l(a_\ga)^{j+1}\,a_\ga^{s\cdot
\al},
$$
Let
$D$ denote the differential operator
$$
D\= (-1)^r\left(\frac \partial{\partial s_1}\dots \frac \partial{\partial s_r}\right).
$$

Let $\hat G(\Ga)$ denote the set of all $\pi\in\hat G$, $\pi\ne triv$ with $N_\Ga(\pi)\ne
0$. For given $\pi\in\hat G$ let $\Lambda(\pi)$ denote the set of all $\la\in\a^*$ with
$m_{\la -\rho}(\pi)\ne 0$. Then $\Lambda(\pi)$ has at most $|W(\h,\g)|$ elements.

Let
$\la\in\a^*$. Since
$\al_1,\dots,\al_r$ is a basis of
$\a^*$ we can write
$\la=\la_1\al_1+\dots+\la_r\al_r$ for uniquely determined $\la_k\in\C$. 

Let $R_k(s)$, $k\in\N$ be a sequence of rational functions on $\C^r$. For an open set
$U\subset\C^r$ let $\N(U)$ be the set of natural numbers $k$ such that the pole-divisor of
$R_k$ does not intersect $U$.
We say that the series
$$
\sum_kR_k(s)
$$
\emph{converges weakly locally uniformly on} $\C^r$ if for every  open $U\subset\C^r$
the series 
$$
\sum_{k\in\N(U)} R_k(s)
$$
converges locally uniformly on $U$.

Let
$$
q_M\df\sum_{p=0}^{\dim\p_M}(-1)^p\dim\(\bigwedge^p\p_M\)^{K_M},
$$
and
$$
q_M^0\df\sum_{p=0}^{\dim\p_M}(-1)^p\dim\(\bigwedge^p\p_M\)^{K_M^0}.
$$

\begin{theorem}\label{2.4}
For $j\in\N$ large enough the series $L^j(s)$ converges locally uniformly in the set
$$
\{ s\in\C : \Re(s_k)>1,\ k=1,\dots, r\}.
$$
The function $L^j(s)$ can be written as Mittag-Leffler series,
\begin{eqnarray*}
L^j(s)&=& D^{j+1} \frac {q_M}{(s_1-1)\cdots (s_r-1)}\\
&&+\sum_{\pi\in\hat G(\Ga)} N_\Ga(\pi)
\sum_{\la\in\Lambda(\pi)} m_{\la-\rho}(\pi) D^{j+1}\frac 1{(s_1+\la_1)\cdots (s_r+\la_r)}.
\end{eqnarray*}
The double series converges weakly locally uniformly on $\C^r$. For $\pi\ne triv$ and
$\la\in\Lambda(\pi)$ we have $\Re(\la_k)> -1$ for $k=1,\dots ,r$.
So in particular, the double series converges
locally uniformly on $\{\Re(s_k)>1\}$.

The same holds for $L^{0,j}(s)$ which satisfies
\begin{eqnarray*}
L^{0,j}(s)&=& D^{j+1} \frac {q_M^0}{(s_1-1)\cdots (s_r-1)}\\
&&+\sum_{\pi\in\hat G(\Ga)} N_\Ga(\pi)
\sum_{\la\in\Lambda(\pi)} m_{\la-\rho}^0(\pi) D^{j+1}\frac 1{(s_1+\la_1)\cdots (s_r+\la_r)}.
\end{eqnarray*}
The integers $q_M, q_M^0$ satisfy
$$
q_M\ \ge\ q_M^0\ >\ 0.
$$
\end{theorem}

The proof will occupy the rest of this section.
We will show that the series $L^j(s)$ converges if the real parts $\Re(s_k)$ are
sufficiently large for $k=1,\dots,r$. Since $L^j(s)$ is a Dirichlet series with positive
coefficients, the convergence in the set $\{\Re(s_k)>1\}$ will follow, once we have
established holomorphy there. This holomorphy will in turn follow from the convergence of
the Mittag-Leffler series.

Since the sum defining $L^{0,j}$ runs over a smaller set, we have for $\Re(s_l)>1$, that $L^j(s)\ge L^{0,j}(s)$.
This implies $q_M\ge q_M^0$.
To see $q_M^0>0$ let $M_\C$ be the complexification of $M$ and let $M_d\subset M_\C$ be a compact form containing $K_M$.
Then $X_{M,d}=M_d/K_M$ is the dual symmetric space to $M/K_M$.
The Betti numbers of $X_{M,d}$ can be computed using the complex $\Omega^\bullet(X_{M,d})^{M_d^0}$ of $M_d^0$-invariant differential forms, where $M_d^0$ is the connected component of the unit.
This complex is isomorphic to
\begin{eqnarray*}
\Omega^\bullet(X_{M,d})^{M_d^0} &\cong & \( C^\infty(M_d^0)\otimes\bigwedge^\bullet \p_M\)^{M_d^0\times K_M^0}\\
&\cong& \left(\bigwedge^\bullet \p_M\)^{K_M^0}.
\end{eqnarray*}
Thus we infer that the Euler number $\chi(X_{M,d})$ equals $q_M^0$.
It is known that Euler numbers of compact symmetric spaces are positive, so $q_M^0>0$.
Note that this deduction of $q_M>0$ is the sole reason for introducing $m_\la^0, q_M^0$, and $L^{0,j}(s)$.

Let
$$
\a_\R^{*,+}\=\{ \la_1\al_1+\dots+\la_r\al_r\mid \la_1,\dots,\la_r >0\}
$$
be the dual positive cone.
Let $\overline{\a_\R^{*,+}}$ be the closure of $\a_\R^{*,+}$ in $\a_\R^*$.

\begin{proposition}\label{2.2}
 Let $\pi\in\hat G$,
$\la\in\a^*$ with
$m_{\la}(\pi)\ne 0$. Then $\Re(\la)$ lies in the set
$$
C\=  -3\rho+\overline{\a_\R^{*,+}}.
$$
For $\pi\in\hat G$ and $\la$ in the boundary of $C$ we have $m_{\la}(\pi)=0$ unless $\pi$ is the
trivial representation and $\la=-3\rho$ in which case $m_\la(\pi)\= q_M$. 
The same assertion holds for $m_\la(\pi)$ replaced with $m_\la^0(\pi)$, only then the integer $q_M$ changes to $q_M^0$.
\end{proposition}

\prf
We introduce a partial order on $\a^*$ by

\begin{tabular}{ccl}
$\mu > \nu$ & $\Leftrightarrow$ & $\mu-\nu$ is a linear combination,\\
&& with positive
integral coefficients, of roots in $\Phi^+$.
\end{tabular}

\begin{lemma}\label{vanishing}
Let $p\in\N$, let $\pi\in\hat G$ and $\mu\in\a^*$ such that $H_p(\n,\pi_K)_\mu\ne 0$. Then
there exists $\nu\in\a^*$ with $\nu<\mu$ and $H_0(\n,\pi_K)_\nu\ne 0$.

Equivalently, if $0\le p<d=\dim(\n)$ and $H^p(\n,\pi_K)_\mu\ne 0$, then there exists
$\eta\in\a^*$ with 
$\eta <\mu$ and
$H^d(\n,\pi_K)_\eta\ne 0$.
\end{lemma}
\prf
The first assertion is a weak version of Proposition 2.32 in \cite{HeSch} and the second
follows from the first and Lemma \ref{1.1}.
\qed

To prove Proposition \ref{2.2} we consider the trivial representation $\pi= triv$ first.
Using the definition of Lie algebra homology it is easy to show that
$m_{-3\rho}(triv)=q_M$ and the other $\la$ with $m_\la(triv)\ne 0$ lie in $-3\rho+\a_R^{*,+}$.
Likewise for $m_\la^0(\pi)$.

For $\pi\ne triv$ we show the stronger statement that if $H^p(\n,\pi_K)_\la\ne 0$, then $\Re(\la)\in -3\rho+{\a_{\R}^{*,+}}$.
We start with the case of $P$ being a minimal parabolic.
Then $M$ is compact, i.e., $M=K_M$ and $q_M=1$.
Using Lemma \ref{1.1} we see that it suffices to show that if $H_0(\n,\pi_K)_\la\ne 0$, then $\Re(\la)\in -\rho+\a_\R^{*,+}$.
So assume $H_0(\n,\pi_K)_\la\ne 0$ and $\pi$ nontrivial.
Theorems 4.16 and 4.25 of \cite{HeSch} imply that $\la$ is a leading coefficient of the asymptotic of matrix coefficients of $\pi$.
By the Howe-Moore Theorem \cite{HM}, these matrix coefficients vanish at infinity on $G$, and this implies that $\Re(\la+\rho)\in\a_\R^{*,+}$.
The case of a minimal parabolic is settled.
 
In general, there is a minimal parabolic $P_0=M_0A_0N_0\subset P=MAN$ such that $M_0\subset M$, $A_0\supset A$, and $N_0\supset N$.
Let $\m_0,\a_0,\n_0$ be the Lie algebras of $M_0$, $A_0$, and $N_0$.
Then
\begin{eqnarray*}
\n_0 &=& \n\oplus \n_M\\
\a_0 &=& \a\oplus \a_M,
\end{eqnarray*}
where $\n_M=\n_0\cap\m$ and $\a_M=\a_0\cap \m$.
Note that $\n$ is an ideal in $\n_0$.

In light of Lemma \ref{vanishing} it suffices to show that if $H^d(\n,\pi_H)\ne 0$, then $\Re(\la)\in -3\rho+\overline{\a_\R^{*,+}}$.
Let $d_M=\dim\n_M$ and $d_0=\dim\n_0$. Then $d_0=d+d_M$.
Consider the Hochschild-Serre spectral sequence
$$
E_2^{p,q}\= H^p(\n_M,H^q(\n,\pi_K))
$$
which abuts to $H^{p+q}(\n_0,\pi_K)$.
We assume that $H^d(\n,\pi_K)_\la\ne 0$.
Then
$$
H^{d_M}(\n_M,H^d(\n,\pi_K)_\la)\ \ne \ 0
$$
as well and thus there exists $\la_M\in\a_M^*$ with
$$
H^{d_M}(\n_M,H^d(\n,\pi_K)_\la)_{\la_M}\ \ne\ 0.
$$
Since $A$ acts trivially on $\n_M$, the latter equals
$$
H^{d_M}(\n_M,H^d(\n,\pi_K))_{\la+\la_M}\= (E_2^{d_M,d})_{\la+\la_M},
$$
where we view $\la+\la_M$ as an element of $\a_0^*=\a^*\oplus\a_M^*$.
The spectral sequence $E$ is supported in the set of indices $0\le p\le d_M$, $0\le q\le d$ and its differentials are $A_0$-homomorphisms.
So $E_2^{d_M,d}$ is the right top corner of this spectral sequence, hence equals $E_\infty^{d_M,d}$ which in this case is $H^{d_0}(\n_0,\pi_K)$.
It follows that $H^{d_0}(\n_0,\pi_K)_{\la+\la_M}\ne 0$ and hence, by the above,.
$$
\Re(\la+\la_M)\ \in\ -3\rho_0 + {\a_{0,\R}^{*,+}},
$$
which by projection implies $\Re(\la)\in -3\rho+{\a_\R^{*,+}}$.
Proposition \ref{2.2} is proved.
\qed 

We continue the proof of Theorem \ref{2.4}. For $a\in A$ set
$$
\ph(a)\= l(a)^{j+1}\ a^{s\cdot \al}.
$$
For $\Re(s_k)>>0$, $k=1,\dots,r$ the Lefschetz formula is valid for this test function. The local
side of the Lefschetz formula equals
$$
\sum_{[\ga]\in\CE_P(\Ga)} \ind(\ga)\, l(a_\ga)^{j+1}\ a_\ga^{s\cdot
\al}\=  L^j(s).
$$
The convergence assertion in the Lefschetz formula implies that the series  converges
absolutely  if
$\Re(s_k)$ is sufficiently large for every
$k=1,\dots r$.  
We will show that it extends to a
holomorphic function in  the set
$\Re(s_k)>1$,
$k=1,\dots ,r$. 
Since $L^j(s)$ is a Dirichlet series with positive coefficients it must
therefore converge in that region.

 With our given test
function and the Haar measure chosen we compute
\eqn
\int_{A^-}\ph(a) a^\la da &=&(-1)^{r(j+1)} \int_{A^-} (\al_1(\log a)\cdots\al_r(\log
a))^{j+1} a^{s\cdot
\al+\la} da\\
\endeqn
\vspace{-30pt}
\eqn
&=& (-1)^{r(j+1)} \int_0^\infty\dots\int_0^\infty \left( {t_1}\cdots t_r\right)^{j+1}
e^{-((s_1+\la_1)t_1+\dots+(s_r+\la_r)t_r)} d t_1\dots d t_r\\
&=& D^{j+1} 
\int_0^\infty\dots\int_0^\infty  e^{-((s_1+\la_1)t_1+\dots+(s_r+\la_r)t_r)} d t_1\dots d t_r\\
&=& D^{j+1}  \frac
1{(s_1+\la_1)\dots(s_r+\la_r)}.
\endeqn
Performing a $\rho$-shift we see that the Lefschetz formula
gives
\begin{eqnarray*}
L^j(s)&=&\sum_{\pi\in\hat G}N_\Ga(\pi)\sum_{\la\in\a^*} m_{\la-\rho}(\pi)\, D^{j+1} \frac
1{(s_1+\la_1)\cdots(s_r+\la_r)}\\
&=& \sum_{\pi\in\hat G}N_\Ga(\pi)\sum_{\la\in\a^*} m_{\la-\rho}(\pi)\,  \frac
{((j+1)!)^r}{(s_1+\la_1)^{j+2}\cdots(s_r+\la_r)^{j+2}}
\end{eqnarray*}
for $\Re(s_k)>>0$.
For every $\pi\in\hat G$ we fix a representative $\La_\pi\in(\a +\t)^*$ of the
infinitesimal character of $\pi$. According to Lemma \ref{1.1}, if
$m_{\la-\rho}(\pi)\ne 0$, then $\la=w\La_\pi|_\a-\rho$ for some $w\in W(\h,\g)$. By abuse of
notation we will write $w\La_\pi$ instead of $w\La_\pi|_\a$. Hence we get
$$
L^j(s)\=\sum_{\pi\in\hat G}N_\Ga(\pi)\sum_{w\in W(\h,\g)} m_{w\La_\pi-2\rho}(\pi)\, D^{j+1}
\frac 1{(s_1+\la_1)\cdots(s_r+\la_r)}.
$$ 

For $\la\in\a^*$ let $\norm\la$ be the norm given by the form $B$ as explained in the
beginning of section \ref{sec1}.

\begin{proposition}\label{2.5}
There are $m\in\N$, $C>0$ such that for every $\pi\in\hat G$ and every $\la\in\a^*$ one has
$$
|m_{\la-\rho}(\pi)|\ \le\ C(1+\norm\la)^m.
$$
\end{proposition}

\prf
Harish-Chandra has shown that there is a locally integrable function $\Theta_\pi^G$ on $G$,
called the global character of $\pi$, such that
$
\tr\pi(h)\= \int_G h(x)\Theta_\pi^G(x)dx
$
for every $h\in C_c^\infty$. It follows that $\Theta_\pi^G$ is invariant under conjugation.
Hecht and Schmid have shown in \cite{HeSch} that for $at\in A^-T$, 
$$
\Theta_\pi^G(at)\= \frac{\sum_{q=0}^{\dim\n} (-1)^q
\Theta_{H_q(\n,\pi_K)}^{AM}(at)}{\det(1-at\mid \n)},
$$
where  $\Theta^{AM}$ is the corresponding global character on the group $AM$.

Let $T=C_1,\dots,C_r$ be a set of representatives of the Cartan subgroups of $M$ modulo $M$-conjugation.
Choose a set of positive roots $\phi_j^+\subset \phi(\c_j,\m)$ for each $j$.
Let $\rho_j=\frac 12\sum_{\al\in\phi_j^+}\al$.
For $x\in C_j$ set
$$
D_{C_j}(x)\= x^{\rho_j}\prod_{\al\in\phi_j^+} (1- x^{-\al}).
$$
This is the Weyl denominator. By the Weyl integration formula the integral
$$
\int_M f_{EP}(m)\sum_{q=0}^{\dim\n}(-1)^{q+\dim\n} \Theta_{H^q(\n,\pi_K)}^{AM}(am)\, dm
$$
equals
$$
 \sum_{j=1}^r\frac {(-1)^{\dim\n}}{|W(C_j,M)|} \int_{C_j^{reg}}D_{C_j}(x) \CO_x(f_{EP})\sum_{q=0}^{\dim\n} (-1)^q\Theta_{H^q(\n\pi_K)}^{AM}(ax)\, dx,
$$
where $f_{EP}$ is the Euler-Poincar\'e function on $M$ and $\CO$ denotes the orbital integral.
Since the orbital integral of the Euler-Poincar\'e function vanishes unless $x$ is elliptic, in which it equal 1 for regular $x$, we see that this equals
$$
\frac {(-1)^{\dim\n}}{|W(T,M)|} \int_{T^{reg}}D_{T}(t) \sum_{q=0}^{\dim\n} (-1)^q\Theta_{H^q(\n\pi_K)}^{AM}(at)\, dt.
$$

On the other hand, by the defining property of the Euler-Poincar\'e function we get that 
$$
\int_M f_{EP}(m)\sum_{q=0}^{\dim\n}(-1)^{q+\dim\n} \Theta_{H^q(\n,\pi_K)}^{AM}(am)\, dm
$$
equals
$$
\sum_{p,q\ge 0} (-1)^{q+\dim\n} \dim\( H^q(\n,\pi_K)\otimes\bigwedge^p\p_M\)^{K_M}
\= \sum_{\al\in\a^*}m_{\la-\rho}(\pi)a^\la.
$$
We put this together  and use the result of Hecht and Schmid to infer
\begin{eqnarray*}
\sum_{\al\in\a^*}m_{\la-\rho}(\pi)a^\la &=& \int_M f_{EP}(m)\sum_{q=0}^{\dim\n}(-1)^{q+\dim\n} \Theta_{H^q(\n,\pi_K)}^{AM}(am)\, dm\\
&=& \frac{(-1)^{\dim\n}}{|W(T,M)|} \int_{T^{reg}} |D_T(t)|^2\sum_{q\ge 0} (-1)^q\Theta_{H^q(\n,\pi_K)}^{AM}(at)\, dt\\
&=& \frac{(-1)^{\dim\n}}{|W(T,M)|}a^{-2\rho} \int_{T^{reg}} |D_T(t)|^2\sum_{q\ge 0} (-1)^q\Theta_{H_q(\n,\pi_K)}^{AM}(at)\, dt\\
&=& \frac{(-1)^{\dim\n}}{|W(T,M)|} a^{-2\rho} \int_{T^{reg}} |D_T(t)|^2 \det(1-at|\n)\Theta_\pi^G(at)\, dt
\end{eqnarray*}

The function $(-1)^{\dim\n}a^{2\rho}\det(1-am|\n)D_T(t)$ equals the Weyl denominator for
$H=AT$. By Theorems 10.35 and 10.48 of \cite{Knapp} there are constants $c_w$, $w\in
W(\h,\g)$ such that
$$
(-1)^{\dim\n}a^{-2\rho}
\det(1-am|\n)D_T(t)\Theta_\pi^G(at)\= \sum_{w\in W(\h\g)}c_w\, (at)^{w\La_\pi}.
$$
We thus have proved the following Lemma.
\begin{lemma}\label{2.6}
 For $a\in A^-$,
$$
\sum_{\la\in\a^*}m_{\la-\rho}(\pi)\, a^\la\= \sum_{w\in
W(\h,\g)}\frac {c_w}{|W(T,M)|}\,a^{w\La_\pi}\int_{T^{reg}}
t^{w\La_\pi-\rho_M}\prod_{\al\in\phi^+(\t,\m)}(1-t^\al)\, dt. 
$$
\end{lemma}

Proposition \ref{2.5} will follow from explicit formulae for the global character
$\Theta_\pi^G$ (see below) which give bounds on the $c_w$.
Another remarkable consequence of Lemma \ref{2.6} is the fact that
there is a finite set $E\subset\t^*$ such that whenever $m_{\la-\rho}(\la)\ne 0$ for some
$\la\in\a^*$ it follows $\La_\pi|_\t\in E$. Hence Proposition \ref{2.5} will follow from the
estimate
$$
|m_{\la-\rho}(\pi)|\ \le\ C(1+\norm{\La_\pi})^m.
$$

In \cite{HC-DSII} Harish-Chandra gives an explicit formula for characters of discrete
series representations which imply the sharper estimate $|m_{\la-\rho}(\pi)|\le C$ for
the discrete series representations. From Harish-Chandra's paper a similar formula can be
deduced for limit of discrete series representations. Alternatively, one can use Zuckerman
tensoring (Prop. 10.44 of \cite{Knapp}) to deduce the estimate for limits of discrete
series representations.
Next, if $\pi=\pi_{\sigma,\nu}$ is induced from some parabolic $P_1=M_1A_1N_1$, then the
character of $\pi$ can be computed from the character of $\sigma$ and $\nu$, see formula
(10.27) in \cite{Knapp}. From this it follows that the claim holds for
standard representations, i.e.
admissible representations which are induced from discrete series or limit of discrete
series representations. 

\begin{lemma}\label{2.7}
There are natural numbers $n,m$ and a constant $d>0$ such that for every $\pi\in\hat G$
there are standard representations $\pi_1,\dots,\pi_n$ and integers $c_1,\dots,c_n$ with
$$
\Theta_\pi\= \sum_{k=1}^n c_k\,\Theta_{\pi_k}
$$
and $|c_k|\le d(1+\norm{\La_\pi}^m)$ for $k=1,\dots, n$.
\end{lemma}

\prf This is Lemma 2.6 of \cite{primgeo}.
\qed

It remains to deduce Theorem \ref{2.4}. Since the coefficients $m_{\la-\rho}(\pi)$ grow at
most like a power of $\norm{\La_\pi}$, the convergence assertion in Theorem \ref{2.4} will
be implied by the following lemma.

\begin{lemma}\label{2.8}
Let $S$ denote the set of all pairs $(\pi,\la)\in\hat G\times\a^*$ such that
$m_{\la-\rho}(\pi)\ne 0$. There is $m_1\in\N$ such that
$$
\sum_{(\pi,\la)\in S} \frac{N_\Ga(\pi)}{(1+\norm{\la})^{m_1}}\ <\ \infty.
$$
\end{lemma}

\prf
By the remark following Lemma \ref{2.6} it suffices to show that 
there is $m\in\N$ such that
$$
\sum_{(\pi,\la)\in S} \frac{N_\Ga(\pi)}{(1+\norm{\La_\pi})^{m_1}}\ <\ \infty.
$$
Let $\pi\in\hat G$. The restriction of $\pi$ to the maximal compact subgroup $K$ decomposes
into finite dimensional isotypes
$$
\pi|_K\=\bigoplus_{\tau\in\hat K} \pi(\tau).
$$
Let $C_K$ be the Casimir operator of $K$ and let
$$
\Delta_G\df -C+2C_K.
$$
Then $\Delta_G$ is the Laplacian on $G$ given by the left invariant metric which at the
point $e\in G$ is given by
$\sp{.}{.}=-B(.,\theta(.))$. Since $\Delta_G$ is left invariant it induces an operator on
$\Ga\bs G$ denoted by the same letter. This operator is
$\ge 0$ and elliptic, so there is a natural number $k$ such that $(1+\Delta_G)^{-k}$ is of
trace class on $L^2(\Ga\bs G)$.
Hence
\begin{eqnarray*}
\infty &>& \tr(1+\Delta_G)^{-k}\\
&=& \sum_{\pi\in\hat G} N_\Ga(\pi) \sum_{\tau\in\hat K}(1-\pi(C)+2\tau(C_K))^{-k}\,
\dim\pi(\tau)\\
&\ge& \sum_{\pi\in\hat G} \frac{N_\Ga(\pi)}{(1- \pi(C)+2\tau_\pi(C_K))^k},
\end{eqnarray*}
where  for each $\pi\in\hat G$ we fix a minimal $K$-type $\tau_\pi$.
Since the infinitesimal character of the minimal $K$-type grows like the infinitesimal
character of $\pi$ the Lemma follows.
\qed

Finally, to prove Theorem \ref{2.4}, let $U\in\C^r$ be open. Let $S(U)$ be the set of all
pairs $(\pi,\la)\in\hat G\times\a^*$ such that $m_{\la-\rho}(\pi)\ne 0$ and the pole divisor
of
$$
\frac 1{(s_1+\la_1)\cdots (s_r+\la_r)}
$$
does not intersect $U$.
Let $V\subset U$ be a compact subset. We have to show that for some $j\in\N$ which does not
depend on $U$ or $V$,
$$
\sup_{s\in V}\sum_{(\pi,\la)\in S(U)} \left| \frac{N_\Ga(\pi) \,
m_{\la-\rho}(\pi)}{(s_1+\la_1)^{j+2}\cdots (s_r+\la_r)^{j+2}}\right|\ <\ \infty.
$$
Let $m$ be as in Lemma \ref{2.5} and $m_1$ as in Lemma \ref{2.8}. Then let $j\ge m+m_1-2$.
Since $V\subset U$ and $V$ is compact there is $\eps>0$ such that $s\in V$ and
$(\pi,\la)\in S(U)$ implies $|s_k+\la_k|\ge\eps$ for every $k=1,\dots,r$. Hence there is
$c>0$ such that for every $s\in V$ and every $(\pi,\la)\in S(U)$,
$$
|(s_1+\la_1)\cdots (s_r+\la_r)|\ \ge\ c(1+\norm\la).
$$
This implies,
\begin{eqnarray*}
\left|\frac{m_{\la-\rho}(\pi)}{(s_1+\la_1)^{j+2}\cdots (s_r+\la_r)^{j+2}}\right|&\le& \frac
1{c^{j+2}}\ \frac{|m_{\la-\rho}(\pi)|}{(1+\norm\la)^{j+2}}\\
&\le& \frac C{c^{j+2}}\ \frac{1}{(1+\norm\la)^{j+2-m}}\\
&\le& \frac C{c^{j+2}}\ \frac{1}{(1+\norm\la)^{m_1}}.
\end{eqnarray*}
The claim now follows from Lemma \ref{2.8}. The proof of Theorem \ref{2.4} is finished.
The version for $L^{0,j}$ is analogous.
\qed

\section{The prime geodesic theorem}\label{sec3}

We now give the two main results of the paper.

\begin{theorem} (Prime Geodesic Theorem)\\
For $T_1,\dots,T_r >0$ let
$$
\Psi(T_1,\dots,T_r)\=\sum_{\stackrel{\stackrel{[\ga]\in\CE_P(\Ga)}{}}{\stackrel{}{a_\ga^{-\al_k}\le
T_k,\ k=1,\dots,r}}}
\la_\ga.
$$
Then, as $T_k\to\infty$ for $k=1,\dots,r$ we have
$$
\Psi(T_1,\dots,T_r)\ \sim\ T_1\cdots T_r.
$$
\end{theorem}
\prf Using Theorem \ref{2.4} the proof is the same as the proof of Theorem 3.1 in \cite{primgeo}.
\qed

Finally, we give a new asymptotic formula for class numbers in number fields. It is quite
different from known results like Siegel's Theorem (\cite{ayoub}, Thm 6.2). The asymptotic is in several
variables and thus contains more information than a single variable one. In a sense it states that
the units of the orders are equally distributed in different directions if only one averages over sufficiently many orders.

Let
$d$ be a prime number
$\ge 3$. Let $r,s\ge 0$ be integers with $d=r+2s$. A number field $F$ is said to be
\emph{of type $(r,s)$} if $F$ has $r$ real and $2s$ complex embeddings.
Let
$S$ be a finite set of primes with
$|S|\ge 2$. Let
$C_{r,s}(S)$ be the set of all number fields $F$ of type $(r,s)$ with the property
$p\in S\
\ \Rightarrow\ \ p$ is non-decomposed in $F$.

Let $O_{r,s}(S)$ denote the set of all orders $\CO$ in number fields $F\in
C_{r,s}(S)$ which are maximal at each $p\in S$. For such an order $\CO$ let $h(\CO)$ be
its class number,
$R(\CO)$ its regulator and $\la_S(\CO)=\prod_{p\in S} f_p$, where $f_p$ is the inertia degree of
$p$ in
$F=\CO\otimes\Q$. Then $f_p\in \{ 1,d\}$ for every $p\in S$.

For $\la\in\CO^\times$ let $\rho_1,\dots,\rho_{r}$ denote the real embeddings of $F$ ordered
in a way that $|\rho_k(\la)|\ge |\rho_{k+1}(\la)|$ holds for $k=1,\dots,r-1$. For the same $\la$
let
$\sigma_1\dots\sigma_{s}$ be pairwise non conjugate complex embeddings ordered in a way
that $|\sigma_k(\la)|\ge |\sigma_{k+1}(\la)|$ holds for $k=1,\dots,s-1$.

For  $k=1,\dots s-1$ let
$$
\al_k(\la)\df 2k(d-2k)\,\log\left( \frac{|\sigma_k(\la)|}{|\sigma_{k+1}(\la)|}\right).
$$
If $s>0$ let 
$$
\al_{s}(\la)\df 2r s \,\log\left( \frac{|\sigma_{s}(\la)|}{|\rho_{1}(\la)|}\right).
$$
For $k=s+1,\dots, r+s-1$ let
$$
\al_k(\la)\df (k+s)(r+s-k)(\,\log\left( \frac{|\rho_{k-s}(\la)|}{|\rho_{k-s+1}(\la)|}\right).
$$

 For $T_1,\dots,T_{r+s-1}>0$ set
$$
v_\CO(T_1,\dots T_{r+s-1})\df \#\{\la\in\CO^\times/\pm 1\mid 0<\al_k(\la)\le T_k,\
k=1,\dots,r+s-1\}.
$$
Let
$$
c=(\sqrt 2)^{1-r-s} \left( \prod_{k=1}^{s-1}
(4k(d-2k)\right)4r s \left(\prod_{k=s+1}^{r+s-1}2(k+s)(r+s-k)\right),
$$
where the factor $4r s$ only occurs if $r s\ne 0$. So $c>0$ and it comes about as
correctional factor between the Haar measure normalization used in the Prime Geodesic
Theorem and the normalization used in the definition of the regulator.

\begin{theorem}
With
$$
\vartheta_S(T)\df \sum_{\CO\in O(S)} v_\CO(T)\, R(\CO)\, h(\CO)\, \la_S(\CO)
$$
we have, as $T_1,\dots,T_{r+s-1}\ra \infty$,
$$
\vartheta(T_1,\dots,T_{r+s-1})\ \sim\ \frac{c}{\sqrt{r+s}}\,T_1\cdots T_{r+s-1}.
$$
\end{theorem}

\prf
For given
$S$ there is a division algebra
$M$ over
$\Q$ of degree
$p$ which splits exactly outside $S$. Fix a maximal order $M(\Z)$ in $M$ and for any ring $R$
define
$M(R)\df M(\Z)\otimes R$. Let $\det :M(R)\ra R$ denote the reduced norm then
$$
\CG(R)\df\{  x\in M(R)\mid \det(x)=1\}
$$
defines a group scheme over $\Z$ with $\CG(\R)\cong\SL_d(\R)=G$. Then $\Ga= \CG(\Z)$ is a
cocompact discrete regular torsion-free subgroup of $G$ (see \cite{class}). As can be seen in
\cite{class}, the theorem can be deduced from the prime geodesic theorem.
\qed

{\small University of Exeter, Mathematics, Exeter EX4
4QE, England\\ a.h.j.deitmar@ex.ac.uk}

\end{document}